\documentclass[10pt]{elsarticle}
\usepackage[cp1251]{inputenc}
\usepackage[english]{babel}
\usepackage{amsmath}
\usepackage{amssymb}
\usepackage{amsfonts}
\usepackage{rotating}
\usepackage{graphicx}
\usepackage{floatflt,epsfig}
\usepackage{lineno,hyperref}
\usepackage{enumerate}
\usepackage{colortbl}
\usepackage[table]{xcolor}
\usepackage{array,tabularx,tabulary,booktabs}
\usepackage{longtable}
\usepackage{multirow}
\usepackage{wrapfig}
\usepackage{subcaption}

\newcolumntype{$}{>{\global\let\currentrowstyle\relax}}
\newcolumntype{^}{>{\currentrowstyle}}

\journal{arXiv}
\newtheorem{thm}{Theorem}
\newtheorem{lem}{Lemma}

\newtheorem{constr}{Construction}

\newcommand{\proof}{\medskip\noindent{\bf Proof.~}}
\begin{document}
\renewcommand{\abstractname}{Abstract}
\renewcommand{\refname}{References}
\renewcommand{\tablename}{Figure.}
\renewcommand{\arraystretch}{0.9}
\thispagestyle{empty}
\sloppy

\begin{frontmatter}
\title{Deza graphs with parameters $(n,k,k-1,a)$ and $\beta = 1$\tnoteref{grant}}
\tnotetext[grant]{
The first, the third and the fourth authors are partially supported by RFBR according to the research project 17-51-560008.
The first and the fourth authors are partially supported by RFBR according to the research project 16-31-00316.
The first author is partially supported by the NSFC (11671258) and STCSM (17690740800).
}

\author[01,03,04]{Sergey~Goryainov}
\ead{44g@mail.ru}

\author[02]{Willem~H.~Haemers\corref{cor1}}
\cortext[cor1]{Corresponding author}
\ead{haemers@uvt.nl}

\author[03]{Vladislav~V.~Kabanov}
\ead{vvk@imm.uran.ru}

\author[03,04]{Leonid~Shalaginov}
\ead{44sh@mail.ru}

\address[01]{Shanghai Jiao Tong University, 800 Dongchuan RD. Minhang District, Shanghai, China}
\address[02]{Dept. of Econometrics and O.R.,
Tilburg University, The Netherlands.}
\address[03]{Krasovskii Institute of Mathematics and Mechanics, S. Kovalevskaja st. 16, Yekaterinburg, 620990, Russia}
\address[04]{Chelyabinsk State University, Brat'ev Kashirinyh st. 129, Chelyabinsk,  454021, Russia}

\begin{abstract}
A Deza graph with parameters $(n,k,b,a)$ is a $k$-regular graph with $n$ vertices in which any two vertices have $a$ or $b$ ($a\leq b$)
common neighbours.
A Deza graph is strictly Deza if it has diameter $2$, and is not strongly regular.
In an earlier paper, the two last authors et el. characterized the strictly Deza graphs with $b=k-1$ and $\beta > 1$,
where $\beta$ is the number of vertices with $b$ common neighbours with a given vertex.
Here we deal with the case $\beta=1$, thus we complete the characterization of strictly Deza graphs with $b=k-1$.
It follows that all Deza graphs with $b=k-1$ and $\beta=1$ can be made from special strongly regular graphs,
and we present several examples of such strongly regular graphs.

A divisible design graph is a special Deza graph, and a Deza graph with $\beta=1$ is a divisible design graph.
The present characterization reveals an error in a paper on divisible design graphs by the second author et al.
We discuss the cause and the consequences of this mistake and give the required errata.
\end{abstract}

\begin{keyword}
Deza graph; divisible design graph; strongly regular graph; dual Seidel switching; involution.
\vspace{\baselineskip}
\MSC[2010] 	05C75\sep 05B30\sep 05E30
\end{keyword}
\end{frontmatter}

\section{Introduction}
A $k$-regular graph $\Gamma$ on $n$ vertices is called a \emph{Deza graph} with parameters $(n,k,b,a)$
if the number of common neighbours of two distinct vertices takes on only two values $a$ or $b$ ($a\leq b$).
If the number of common neighbours of two vertices only depend on whether the vertices are adjacent or not,
then $\Gamma$ is a strongly regular graph with parameters $(n,k,\lambda,\mu)$, where $\lambda$ ($\mu$) is
the number of common neighbours of two adjacent (non-adjacent) vertices; so $\{a,b\} = \{\lambda,\mu\}$.
A Deza graph is called a \emph{strictly Deza graph} if it has diameter $2$ and is not strongly regular.
Note that the complete graph $K_n$ (which is normally excluded from being strongly regular) is a Deza
graph which is not strictly Deza because it has diameter~$1$.

Let $\Gamma$ be a Deza graph with parameters $(n, k, b, a)$, and let $v$ be a vertex of $\Gamma$.
Denote by $N(v)$ the set of neighbours of a vertex $v$, and let $\beta(v)$ be the number of vertices $u\in V(\Gamma)$
such that $|N(v)\cap N(u)| = b$.
\begin{lem}[\cite{EFHHH98}, Proposition~1.1]\label{beta} The number $\beta(v)$ does not depend on the choice of $v$ and is given by
$$\beta(v) = \beta = \frac{k(k - 1) - a(n - 1)}{b - a}\ \mbox{if}\ a\neq b,\ \mbox{and}\ \beta = n-1\ \mbox{if}\ a=b.$$
\end{lem}

Strictly Deza graphs with parameters $(n, k, b, a)$, where $k = b + 1$ and $\beta > 1$ hold, were investigated in \cite{KMSh}.
The following theorem was proved.
\begin{thm}[\cite{KMSh}]\label{>1}
Let $\Gamma$ be a strictly Deza graph with parameters $(n, k, b, a)$ and  $\beta > 1$.
The parameters $k$ and $b$ of $\Gamma$ satisfy the condition $k = b + 1$  if and only if $\Gamma$ is isomorphic to the strong product
of $K_2$ with the complete multipartite graph with parts of size $(n-k+1)/2$. 
\end{thm}

In this paper we characterize strictly Deza graphs with parameters $(n, k, b, a)$, where $k = b + 1$ and $\beta = 1$.

Note that the adverb \lq{strictly{\rq} in Theorem~\ref{>1} can not be removed, as is shown by the $n$-cycle with $n\geq 5$.
However, in the present characterization, we will see that there are no examples which are not strictly Deza, except for $K_2$.

\section{The characterization}

We present two constructions of Deza graphs with parameters $(n, k, b, a)$, where $k = b + 1$ and $\beta = 1$.
Both constructions use a strongly regular graph $\Delta$ with parameters $(m,\ell,\lambda,\mu)$ where $\lambda = \mu-1$.

\begin{constr}\label{C1}
Let $\Gamma_1$ be the strong product of $K_2$ and $\Delta$.
The graph $\Gamma_1$ is a strictly Deza graph with parameters $(n,k,k-1,a)$ and $\beta = 1$, where $n=2\ell,~k=2\ell+1,~a=2\mu$.
\end{constr}

So, if $B$ and $A_1$ are the adjacency matrices of $\Delta$ and $\Gamma_1$, respectively,
then $A_1=B\otimes J_2 - I_{n}$ ($J_m$ is the $m\times m$ all-ones matrix, and $I_m$ is the identity matrix of order $m$).

Suppose that $\Delta$ has an involution that interchanges only non-adjacent vertices.
Let $P$ be the corresponding permutation matrix, then $B'=PB$ is a symmetric matrix (because $P=P^\top$ and $PBP=B$)
with zero diagonal (because $P$ interchanges only nonadjacent vertices).
So $B'$ is the adjacency matrix of a graph $\Delta'$ (say), which is a Deza graph because $B'^2=PBPB=B^2$.
This construction was given in \cite{EFHHH98} and the method has been called {\em dual Seidel switching}; see~\cite{H84}.

Next, let $\Gamma_2'$ be the strong product of $K_2$ and $\Delta'$.
Modify $\Gamma_2'$ as follows:
For any transposition $(x~y)$ of the involution, take the corresponding two pairs of vertices $x',x''$ and $y',y''$ in $\Gamma_2'$,
delete the edges $\{x',x''\}$ and $\{y',y''\}$, and insert the edges $\{x',y''\}$ and $\{x'',y'\}$.
Define $\Gamma_2$ to be the resulting graph.
If $A_2$ is the adjacency matrix of $\Gamma_2$, then we can also construct $\Gamma_2$ from $\Gamma_1$ using dual Seidel switching
in the following way: $A_2=P_1A_1$ where $P_1=P\otimes I_2$.
We easily have that $A_2^2=A_1^2$, which shows that $\Gamma_2$ is a Deza graph with the same parameters as $\Gamma_1$.

\begin{constr}\label{C2}
The graph $\Gamma_2$ is a strictly Deza graph with parameters $(n,k,k-1,a)$ and $\beta = 1$, where $n=2m,~k=2\ell+1,~a=2\mu$.
\end{constr}

Note that in $\Gamma_1$ any two vertices with $b$ common neighbors are adjacent.
For $\Gamma_2$ this is not true, therefore $\Gamma_1$ and $\Gamma_2$ are non-isomorphic.

\begin{thm}\label{mainthm} If $\Gamma$ is a Deza graph with parameters $(n,k,k-1,a)$, $k > 1$, and $\beta = 1$,
then $\Gamma$ can be obtained either from Construction~\ref{C1} or from Construction~\ref{C2}.
\end{thm}

In case $k=1$, $\Gamma$ consists of $n/2$ disjoint edges and $\beta=1$ implies $\Gamma=K_2$.

\section{Proof of the characterization}

Let $\Gamma$ be a Deza graph with parameters $(n,k,k-1,a)$ with $k>1$ and $\beta = 1$.
For a vertex $x$ of $\Gamma$, denote by $x_b$ the vertex of $\Gamma$ that has $b$ common neighbours with $x$.
Note that $(x_b)_b = x$ holds.
A vertex $x$ in $\Gamma$ is said to be an $A$\emph{-vertex} ($NA$\emph{-vertex}) if $x$ is adjacent (not adjacent) to the
vertex $x_b$.

\begin{lem}\label{ComB}
Let $x$ be an $A$-vertex. Then the equality $N(x) \setminus \{x_b\} = N(x_b) \setminus \{x\}$ holds.
\end{lem}

\begin{lem}\label{TwoAPairs} Let $x,y$ be two $A$-vertices, $\{x,x_b\} \ne \{y, y_b\}$. Then
there are either all possible edges between $\{x,x_b\}$ and $\{y, y_b\}$ or no such edges.
\end{lem}
\proof It follows from Lemma \ref{ComB}. $\square$

\begin{lem}\label{APaierNAPair} Let $x$ be an $A$-vertex and $y$ be a $NA$-vertex in $\Gamma$.
Then there are either all possible edges between $\{x,x_b\}$ and $\{y, y_b\}$ or no such edges.
\end{lem}
\proof Suppose $x$ is adjacent to $y$. Since $x$ is an $A$-vertex, $x_b$ is adjacent to $y$.
Since the vertices $y,y_b$ have $k-1$ common neighbours and $y$ is adjacent to $x$ and $x_b$,
the vertex $y_b$ is adjacent to at least one of the vertices $x$ and $x_b$.
The fact that $x,x_b$ are $A$-vertices implies that $y_b$ is adjacent to both the vertices $x$ and $x_b$. $\square$

\medskip
Let $x$ be a $NA$-vertex. Then $N(x)$ contains precisely one vertex which is not adjacent to $x_b$.
Denote this vertex by $x'$.

\begin{lem} Let $x$ be a $NA$-vertex. Then the vertex $(x_b)'$ belongs to $N_2(x)$.
\end{lem}
\proof By definition, $(x_b)'$ is the neighbour of $x_b$ which is not adjacent to $(x_b)_b = x$. $\square$

\begin{lem}\label{W}
The following statements hold.\\
{\rm(1)} The vertex $(x_b)'$ is a unique neighbour of $x_b$ in $N_2(x)$.\\
{\rm(2)} The common neighbours of $x_b$ and $(x_b)'$ lie in $N(x)$.\\
{\rm(3)} The vertex $(x_b)'$ has precisely $a$ neighbours in $N(x)$.\\
{\rm(4)} Any vertex $u \in N(x,x')$ is adjacent to $(x_b)'$.\\
{\rm(5)} The equality $N(x,x') = N(x_b,x') = N(x, (x_b)') = N(x',(x_b)')$ holds.\\
{\rm(6)} If a vertex is adjacent to any three ones from the set $\{x,x',x_b,(x_b)'\}$,
the is adjacent to all of them.
\end{lem}
\proof
(1) It follows from the fact that the vertex $(x_b)'$ has precisely $b = k-1$ neighbours in $N(x)$.\\
(2) It follows from item (1).\\
(3) These $a$ neighbours are the vertices from $N(x,(x_b)')$.\\
(4) The vertex $u$ has precisely $a$ neighbours in $N(x)$; one of them is $x'$.
The vertex $x_b$ is adjacent to all vertices in $N(x)$ but $x'$. This means
that $u$ and $x_b$ has precisely $a-1$ common neighbors in $N(x)$.
The only neighbour of $x_b$ in $N_2(x)$ is the vertex $(x_b)'$.
Thus, $(x_b)'$ is a common neighbour of $u$ and $x_b$, and, in particular, $u$ is adjacent to $(x_b)'$.\\
(5) It follows from the fact that any vertex from $N(x,x')$ is adjacent to $x_b$ and $(x_b)'$.\\
(6) It follows from item (5).
$\square$

\medskip
For a $NA$-vertex $x$, put $W(x) = N(x,x')$.

\begin{lem}\label{W1} For any $NA$-vertex $x$ and for any vertex $y \in W(x)$,
the vertex $y_b$ belongs to $W(x)$.
\end{lem}
\proof
Since the vertices $x,x',x_b,(x_b)'$ are neighbours of $y$, and $|N(y,y_b)| = k-1$,
the vertex $y_b$ is adjacent to at least three of them. By Lemma \ref{W}(6),
the vertex $y_b$ is adjacent to all of them. $\square$

\begin{lem}\label{aiseven}
The parameter $a$ is even.
\end{lem}
\proof  It follows from Lemma \ref{W1} and the fact that $|W(x)| = a$. $\square$

\begin{lem}
For any $NA$-vertex $x$, the vertex $x'$ is a $NA$-vertex.
\end{lem}
\proof Suppose that $(x')_b$ and $x'$ are adjacent. Since $x$ is a neighbour of $x'$,
then, in view of Lemmas \ref{TwoAPairs} and \ref{APaierNAPair}, the vertices $x,x_b$ belong to $N(x',(x')_b)$.
In particular, we obtain that $x'$ is adjacent to $x_b$, which is a contradiction because $x'$
is not adjacent to $x_b$ by definition. $\square$

\begin{lem}\label{xpb}
For a $NA$-vertex $x$, the vertex $(x')_b$ is not adjacent to $x$.
\end{lem}
\proof Suppose that $(x')_b$ belongs to $N(x)$. Since $(x')_b \ne x'$, the vertex $(x')_b$
is adjacent to $x_b$. This means that $x_b$ is a unique vertex which is adjacent to $(x')_b$ and is not adjacent to $x'$.
Thus, every vertex $y$ from $N(x,(x')_b)$ is adjacent to $x'$, and, consequently, belongs to $N(x,x') = W(x)$.
Since $|W(x)| = |N(x,(x')_b)| = a$, we have $W(x) = N(x,(x')_b)$. In particular, this gives the inclusion
$W(x) \subseteq N((x')_b)$.

Let us count the number of common neighbours of the vertices $(x')_b$ and $(x_b)'$. We have
$a = |N((x')_b, (x_b)')| \ge |W(x) \cup \{x_b\}| = a+1$, which is a contradiction. $\square$

\begin{lem}\label{xppeqx}
For a $NA$-vertex $x$, the equality $x'' = x$ holds.
\end{lem}
\proof The vertex $x'$ is adjacent to $x$ be definition. By Lemma \ref{xpb}, the vertex $(x')_b$ is not adjacent to $x$.
This proves the lemma. $\square$

\begin{lem}\label{xpbeqxbp}
For a $NA$-vertex $x$, the equality $(x')_b = (x_b)'$ holds.
\end{lem}
\proof
It is enough to show that $(x')_b$ and $(x_b)'$ have at least $a+1$ common neighbours.

By Lemma \ref{W}(5), the vertex $(x_b)'$ is adjacent to all vertices in $W(x) = N(x,x')$.
By definition of $(x')_b$, in view of Lemma \ref{xppeqx}, the vertex $(x')_b$ is adjacent to all vertices in $W(x) = N(x,x')$.
Thus, the vertices $(x')_b$ and $(x_b)'$ have at least $|W(x)|=a$ common neighbours.

The vertex $x'$ has $k-a-1$ neighbours in $N_2(x)$, which are neighbours of $(x')_b$ too.
Let us take a vertex $y\in N(x',(x')_b) \cap N_2(x)$. Since $y$ is adjacent to $x'$ and
$x_b$ is not adjacent to $x'$, the vertices $y$ and $x_b$ have precisely
$a-1$ common neighbours in $N(x)$. This implies that the vertex $(x_b)'$,
which is a unique neighbour of $x_b$ in $N_2(x)$, is a common neighbour of $x_b$ and $y$,
and, in particular, $y$ is adjacent to $(x_b)'$. Thus, the vertices $(x')_b$ and $(x_b)'$
have at least $a+1$ common neighbours, which proves the lemma. $\square$

\medskip
Further, in view of Lemma \ref{xpbeqxbp}, we use the simplified notation $x_b'$ for the vertex $(x')_b = (x_b)'$.
For a $NA$-vertex $x$, put $C(x) = \{x,x',x_b,x_b'\}$.

\begin{lem}\label{C}
For a $NA$-vertex $x$, the equalities $C(x) = C(x') = C(x_b) = C(x_b')$ hold.
\end{lem}
\proof It follows from the equality $(x_b)_b = x$, Lemma \ref{xppeqx} and Lemma \ref{xpbeqxbp}.

\begin{lem}\label{quadr}
The set $NA$-vertices in $\Gamma$ can be partitioned into quadruples of the form $\{x,x',x_b,x_b'\}$.
\end{lem}
\proof It follows from Lemma \ref{C}. $\square$

\begin{lem}\label{xxbyyb} Let $x$ be a $NA$-vertex, and let $y$ be a vertex, $y$ does not belong to $\{x,x',x_b,x_b'\}$.
Then there are either all possible edges between $\{x,x_b\}$ and $\{y,y_b\}$ or no such edges.
\end{lem}
\proof If $y$ is an $A$-vertex, then the result follows from Lemma \ref{APaierNAPair}. Let us assume without loosing
of generality that $y$ is a $NA$-vertex. In view of Lemma \ref{C},
it is enough to show that if $x$ and $y$ are adjacent, then $x$ is adjacent to $y_b$,
$x_b$ is adjacent to $y$, and $x_b$ is adjacent to $y_b$.

Suppose $x$ is not adjacent to $y_b$. Then $x = y'$ holds, which is a contradiction since,
by Lemma \ref{quadr}, the condition $\{x,x',x_b,x_b'\} \cap \{y,y',y_b,y_b'\} = \emptyset$ holds.
Using the similar arguments, we can show that $x_b$ is adjacent to $y$ and $x_b$ is adjacent to $y_b$.
The lemma is proved. $\square$

\begin{lem}\label{xpyp}
For any two $NA$-vertices $x,y$ such that $C(x) \ne C(y)$,
the following statements hold.\\
{\rm (1)} The vertices $x'$ and $y$ are adjacent if and only if the vertices $x$ and $y'$ are adjacent.\\
{\rm (2)} The vertices $x$ and $y$ are adjacent if and only if the vertices $x'$ and $y'$ are adjacent.
\end{lem}
\proof (1) Suppose there exists two $NA$-vertices $x,y$, $\{x,x',x_b,x_b'\} \ne \{y,y',y_b,y_b'\}$
such that $x'$ is adjacent to $y$ and $x$ is not adjacent $y'$. For any vertex $z \in N(x,y)$, $z \ne x',y'$,
in view of Lemma \ref{xxbyyb}, there are all possible edges between $\{z,z_b\}$ and $\{x,x_b\}$, $\{z,z_b\}$ and $\{y,y_b\}$.
In particular, $z_b$ belongs to $N(x,y)$. Moreover, $x'$ belongs to $N(x,y)$ and $y'$ does not belong to $N(x,y)$.
This means that $a = |N(x,y)|$ is odd, which contradicts to Lemma \ref{aiseven}.
 \\
(2) It follows from item (1). $\square$

\begin{lem}\label{CxCy}
Let $x,y$ be two $NA$-vertices, $C(x) \ne C(y)$. Then the following statements hold.\\
{\rm(1)} If a vertex from $\{x,x_b\}$ is adjacent to a vertex from $\{y,y_b\}$ or $\{x',x_b'\}$ is adjacent to a vertex from $\{y',y_b'\}$,
then there are all possible edges between $\{x,x_b\}$ and $\{y,y_b\}$,
and there are all possible edges between $\{x',x_b'\}$ and $\{y',y_b'\}$.\\
{\rm(2)} If a vertex from $\{x,x_b\}$ is adjacent to a vertex from $\{y',y_b'\}$
or a vertex from $\{x',x_b'\}$ is adjacent to a vertex from $\{y,y_b\}$,
then there are all possible edges between $\{x,x_b\}$ and $\{y',y_b'\}$,
and there are all possible edges between $\{x',x_b'\}$ and $\{y,y_b\}$.\\
{\rm(3)} There are either $0$, or $8$, or $16$ edges between the sets $C(x)$ and $C(y)$.
\end{lem}
\proof (1) Without loosing of generality, assume that $x$ and $y$ are adjacent.
By Lemma \ref{xxbyyb}, there are all possible edges between $\{x,x_b\}$ and $\{y,y_b\}$.
On the other hand, by Lemma \ref{xpyp}(2), the vertices $x'$ and $y'$ are adjacent.
Then, by Lemma \ref{xxbyyb} again, there are all possible edges between $\{x',x_b'\}$ and $\{y',y_b'\}$.\\
(2) The proof is similar to the one from item (1).\\
(3) It follows from items (2) and (3). $\square$

\begin{lem}\label{NA_A}
Let $x$ be a $NA$-vertex and $y$ be an $A$-vertex. The following statements hold.\\
{\rm(1)} If $x$ and $y$ are adjacent, then the vertices $y$ and $y_b$ lie in $W(x)$.\\
{\rm(2)} There are either all possible edges between the sets $\{x,x',x_b,x_b'\}$ and $\{y,y_b\}$,
or no such edges.
\end{lem}
\proof
(1) Let $z$ be a vertex in $N(x,y)$, $z\ne x',y_b$. We prove that $z_b$ belongs to $N(x,y)$.
If $z$ is an $A$-vertex, then it follows from Lemma \ref{APaierNAPair}.
Suppose $z$ is an $NA$-vertex. Then, by Lemma \ref{APaierNAPair}, $z_b$ belongs to $N(y)$, and
by Lemma \ref{xxbyyb}, $z_b$ belongs to $N(x)$, which implies that $z_b$ belongs to $N(x,y)$.
Moreover, $y_b$ belongs to $N(x,y)$. In view of Lemma \ref{aiseven}, $|N(x,y)|$ is even,
which implies that $x'$ belongs to $N(x,y)$, and, in particular, $y$ is adjacent to $x'$.
This means that $y$ belongs to $N(x,x')$, and, consequently, $y_b$ belongs to $N(x,x')$.\\
(2) It follows from item (1). $\square$

\medskip
Let $\Gamma'$ be the graph obtained from $\Gamma$ by removing all edges $\{x,x_b\}$, where $x$ is an $A$-vertex,
and all edges $\{x,x'\}$, where $x$ is a $NA$-vertex.

\begin{lem}\label{Gammap}
The following statements hold.\\
{\rm(1)} The graph $\Gamma'$ is $(k-1)$-regular.\\
{\rm(2)} For any vertices $x,y$ in $\Gamma'$ such that $y \ne x_b$,
there are either all possible edges between $\{x,x_b\}$ and $\{y,y_b\}$ or no such edges.\\
{\rm(3)} For any vertices $x,y$ in $\Gamma'$ such that $y \ne x_b$,
the equalities $N_{\Gamma'}(x,y) = N_{\Gamma'}(x_b,y) = N_{\Gamma'}(x,y_b) = N_{\Gamma'}(x_b,y_b)$ hold.
\end{lem}
\proof
(1) It follows from the fact that precisely one edge was removed for each vertex.\\
(2) If $x, y$ are $A$-vertices, then it follows from Lemma \ref{TwoAPairs}.
If one of the vertices $x,y$ is an $A$-vertex and the other is a $NA$-vertex,
then it follows from Lemma \ref{NA_A}. If $x, y$ are $NA$-vertices, then it follows from Lemma \ref{CxCy}.\\
(3) It follows from the fact that in the graph $\Gamma'$ any vertex $x$ has the same neighbourhood as the vertex $x_b$.
$\square$

\medskip
Let $\Gamma''$ be the graph whose vertex set is the set of all pairs of vertices $\{x,x_b\}$ in $\Gamma'$,
and two vertices $\{x,x_b\}$ and $\{y,y_b\}$ are adjacent in $\Gamma''$
whenever there are all possible edges between then sets of vertices $\{x,x_b\}$ and $\{y,y_b\}$ in $\Gamma'$.
\begin{lem}\label{GammappIsDeza} The graph $\Gamma''$ is a Deza graph with parameters
$(\frac{v}{2},\frac{k-1}{2},\frac{a}{2},\frac{a-2}{2})$,
where any two adjacent (distinct non-adjacent) vertices $\{x, x_b\}$ and $\{y, y_b\}$ have $\frac{a}{2}$ ($\frac{a-2}{2}$, respectively)
common neighbours iff the vertices $x$ and $y$ are $NA$-vertices, $\{x', x_b'\} \ne \{y, y_b\}$ holds, and there are 8 edges between
the sets $\{x,x',x_b,x_b'\}$ and $\{y, y', y_b, y_b'\}$.
\end{lem}
\proof It follows from Lemma \ref{Gammap}(1) that the graph $\Gamma''$ is $(k-1)/2$-regular.
For arbitrary two distinct vertices $\{x, x_b\}$ and $\{y, y_b\}$ in $\Gamma''$,
we consider all possible cases and prove that the number of their common neighbours is either $a/2$ or $(a-2)/2$.

(1) The vertices $x$ and $y$ are $A$-vertices.
By Lemma \ref{Gammap}(2),
there are either all possible edges between $\{x,x_b\}$ and $\{y, y_b\}$ or no such edges.

Suppose there are all possible edges between $\{x,x_b\}$ and $\{y, y_b\}$.
Then $N_\Gamma(x,y) = \{x_b, y_b\} \cup N_{\Gamma'}(x,y)$, which, in view of Lemma \ref{Gammap}(3), implies
that the adjacent vertices $\{x, x_b\}$ and $\{y, y_b\}$ have $\frac{|N_{\Gamma'}(x,y)|}{2} = \frac{a-2}{2}$ common neighbours
in $\Gamma''$.

Suppose there are no edges between $\{x,x_b\}$ and $\{y, y_b\}$.
Then $N_\Gamma(x,y) = N_{\Gamma'}(x,y)$, which, in view of Lemma \ref{Gammap}(3), implies
that the non-adjacent vertices $\{x, x_b\}$ and $\{y, y_b\}$ have $\frac{|N_{\Gamma'}(x,y)|}{2} = \frac{a}{2}$ common neighbours
in $\Gamma''$.

(2) The vertex $x$ is an $NA$-vertex and the vertex $y$ is an $A$-vertex. By Lemma \ref{NA_A}(2),
there are either all possible edges between $\{x,x',x_b,x_b'\}$ and $\{y, y_b\}$ or no such edges.

Suppose there are all possible edges between $\{x,x',x_b,x_b'\}$ and $\{y, y_b\}$.
Then $N_\Gamma(x,y) = \{x', y_b\} \cup N_{\Gamma'}(x,y)$, which, in view of Lemma \ref{Gammap}(3), implies
that the adjacent vertices $\{x, x_b\}$ and $\{y, y_b\}$ have $\frac{|N_{\Gamma'}(x,y)|}{2} = \frac{a-2}{2}$ common neighbours
in $\Gamma''$.

Suppose there are no edges between $\{x,x',x_b,x_b'\}$ and $\{y, y_b\}$.
Then $N_\Gamma(x,y) = N_{\Gamma'}(x,y)$, which, in view of Lemma \ref{Gammap}(3), implies
that the non-adjacent vertices $\{x, x_b\}$ and $\{y, y_b\}$ have $\frac{|N_{\Gamma'}(x,y)|}{2} = \frac{a}{2}$ common neighbours
in $\Gamma''$.

(3) The vertices $x$ and $y$ are $NA$-vertices. If $C(x) = C(y)$ (in other words, $\{y,y_b\} = \{x', x_b'\}$),
then $N_\Gamma(x,y) = N_{\Gamma'}(x,y)$, which, in view of Lemma \ref{Gammap}(3), implies
that the non-adjacent vertices $\{x, x_b\}$ and $\{y, y_b\}$ have $\frac{|N_{\Gamma'}(x,y)|}{2} = \frac{a}{2}$ common neighbours
in $\Gamma''$. If $C(x) \ne C(y)$ (in other words, $\{y,y_b\} \ne \{x', x_b'\}$), then,
by Lemma \ref{CxCy}, there are either $0$, or $8$, or $16$ edges between the sets $C(x) = \{x,x',x_b,x_b'\}$ and
$C(y) = \{y,y',y_b,y_b'\}$. Let us consider these cases.

Suppose there are no edges between $\{x,x',x_b,x_b'\}$ and $\{y,y',y_b,y_b'\}$.
Then $N_\Gamma(x,y) = N_{\Gamma'}(x,y)$, which, in view of Lemma \ref{Gammap}(3), implies
that the non-adjacent vertices $\{x, x_b\}$ and $\{y, y_b\}$ have $\frac{|N_{\Gamma'}(x,y)|}{2} = \frac{a}{2}$ common neighbours
in $\Gamma''$.

Suppose there are all 16 possible  edges between $\{x,x',x_b,x_b'\}$ and $\{y,y',y_b,y_b'\}$.
Then $N_\Gamma(x,y) = \{x', y'\} \cup N_{\Gamma'}(x,y)$, which, in view of Lemma \ref{Gammap}(3), implies
that the adjacent vertices $\{x, x_b\}$ and $\{y, y_b\}$ have $\frac{|N_{\Gamma'}(x,y)|}{2} = \frac{a-2}{2}$ common neighbours
in $\Gamma''$.

Suppose there are 8 edges between $\{x,x',x_b,x_b'\}$ and $\{y,y',y_b,y_b'\}$. Then there are two subcases with respect
to the items (1) and (2) of Theorem \ref{CxCy}. If there are all possible edges between $\{x,x_b\}$ and $\{y,y_b\}$
and there are all possible edges between $\{x',x_b'\}$ and $\{y',y_b'\}$,
then $N_\Gamma(x,y) = N_{\Gamma'}(x,y)$, which, in view of Lemma \ref{Gammap}(3), implies
that the adjacent vertices $\{x, x_b\}$ and $\{y, y_b\}$ have $\frac{|N_{\Gamma'}(x,y)|}{2} = \frac{a}{2}$ common neighbours
in $\Gamma''$. If there are all possible edges between $\{x,x_b\}$ and $\{y',y_b'\}$
and there are all possible edges between $\{x',x_b'\}$ and $\{y,y_b\}$,
then $N_\Gamma(x,y) = \{x', y'\} \cup N_{\Gamma'}(x,y)$, which, in view of Lemma \ref{Gammap}(3), implies
that the non-adjacent vertices $\{x, x_b\}$ and $\{y, y_b\}$ have $\frac{|N_{\Gamma'}(x,y)|}{2} = \frac{a-2}{2}$ common neighbours
in $\Gamma''$.
$\square$

\medskip
By Lemma \ref{CxCy}(3), for any two $NA$-vertices $x,y$ in $\Gamma$ such that $C(x) \ne C(y)$,
there are either $0$, or $8$, or $16$ edges between the sets $C(x) = \{x,x',x_b,x_b'\}$ and $C(y) = \{y,y',y_b,y_b'\}$.
This implies that there are either 0, or 2, or $4$ edges between the sets of vertices $\{\{x,x_b\},\{x',x_b'\}\}$
and $\{\{y,y_b\},\{y',y_b'\}\}$. Note that, in the case when there are $2$ edges,
the switching edges between $\{\{x,x_b\},\{x',x_b'\}\}$
and $\{\{y,y_b\},\{y',y_b'\}\}$ preserves regularity of $\Gamma''$. Let us make all
such switchings in $\Gamma''$. Denote by $\Gamma'''$ the resulting graph.

\begin{lem}\label{Gammappp} The following statements hold.\\
{\rm(1)} The adjacency matrix of $\Gamma'''$ can be obtained from the adjacency matrix of $\Gamma''$
by permuting rows in pairs corresponding to the vertices $\{x,x_b\},\{x',x_b'\}$
for all quadruples $\{x,x_b, x',x_b'\}$ of $NA$-vertices.
In particular, the adjacency matrices of $\Gamma'''$ and $\Gamma''$ coincide in the case when
$\Gamma$ has no $NA$-vertices.\\
{\rm(2)} The graph $\Gamma'''$ is strongly regular with parameters $(\frac{v}{2},\frac{k-1}{2},\frac{a-2}{2},\frac{a}{2})$.
\end{lem}
\proof (1) It follows from the definition of $\Gamma'''$ and Lemmas \ref{CxCy} and \ref{NA_A}(2).\\
(2) It follows from Lemma \ref{GammappIsDeza} and item (1). $\square$

\medskip
In view of Lemma \ref{Gammappp}, if $\Gamma$ has no $NA$-vertices, then it comes from Construction \ref{C1}.
Let us prove that, if $\Gamma$ has $NA$-vertices, then it comes from Construction \ref{C2}.
Lemmas \ref{GammappIsDeza} and \ref{Gammappp} imply that the adjacency matrix of the Deza graph $\Gamma''$ can be obtained from
the adjacency
matrix of the strongly regular $\Gamma'''$ by swapping rows in pairs corresponding to the vertices
$\{x,x_b\}$ and $\{x',x_b'\}$ for all quadruples $\{x,x',x_b, x_b'\}$ of $NA$-vertices.
It follows from \cite[Theorem 3.1]{EFHHH98}, that this permutation is an order 2  automorphism of $\Gamma'''$
that interchanges only non-adjacent vertices. Thus, $\Gamma''$ can be obtained from $\Gamma'''$
by dual Seidel switching, which implies that $\Gamma$ comes from Construction \ref{C2} by definition.
The theorem is proved.

\section{Strongly regular graphs with $\lambda=\mu-1$}

Note that, if $\Gamma$ is a strongly regular graph with $\lambda=\mu-1$, then so is its complement.
So both $\Gamma$ and its complement satisfy the condition for Construction~\ref{C1}.
For Construction~\ref{C2}, $\Gamma$ needs an involution that interchanges only nonadjacent vertices.
In this section we survey strongly regular graphs with $\lambda=\mu-1$, and look for the desired involutive
automorphism.

\subsection{Paley graphs of square order}

The Paley graph $P(r)$ is a graph with vertex set $\mathbb{F}_r$, where $r$ is a prime power such that $r\equiv 1\mod 4$.
Two vertices $x$ and $y$ of $P(r)$ are adjacent whenever $x-y$ is a non-zero square in $\mathbb{F}_r$.
See~\cite{J17} for an excellent survey of Paley graphs.
The Paley graph is a strongly regular graph with parameters $(r, \frac{r-1}{2}, \frac{r-5}{4}, \frac{r-1}{4})$,
so it satisfies the conditions of Construction~\ref{C1}, which leads to strictly Deza graphs with parameters
$(2r,r,r-1,(r-1)/2$).
The complement of $P(r)$ is isomorphic to $P(r)$, since for any non-square $a$ in $\mathbb{F}_r$ the map $x\rightarrow ax$
interchanges edges and non-edges in $P(r)$.

If $r=q^2$ is a square, then the Paley graph $P(q^2)$ satisfies the conditions of Construction~\ref{C2}.
To explain this we need some properties of the field $\mathbb{F}_{q^2}$.
Let $d$ be a non-square in $\mathbb{F}_{q}^*$. The elements of the finite field of order $q^2$ can be considered as
$$\mathbb{F}_{q^2} =
\{x+y\alpha~|~x,y \in \mathbb{F}_q\},$$ where $\alpha$ is a root of the polynomial $f(t) = t^2 - d.$

Let $\beta$ be a primitive element of the finite field $\mathbb{F}_{q^2}.$
Then we have $\mathbb{F}_{q}^* = \{\beta^{i(q+1)}~|~ i \in \{0, \ldots,q-2\}\}$.
Since $q+1$ is even, each element  of $\mathbb{F}_{q}^*$ is a square in $\mathbb{F}_{q^2}^*$.
It also follows that $x^{q-1}=\beta^{q^2-1}=1$  for every $x\in\mathbb{F}_q^*$.

\begin{lem}\label{Frobenius}
For any $\gamma = x + y\alpha$ from $\mathbb{F}_{q^2}$, the following equalities hold.\\
{\rm(1)} $\gamma^q = x - y\alpha$;\\
{\rm(2)} $ \gamma - \gamma^q = 2y\alpha$.
\end{lem}
\proof
(1)$(x+y\alpha)^q=x^q+y^q\alpha^q=x-y\alpha$.
(2) It follows from item (1). $\square$\\

For any $\gamma = x+y\alpha\in \mathbb{F}_{q^2}^*$ define the \emph{norm} mapping $N$ by $N(\gamma) = \gamma^{q+1} =
\gamma\gamma^{q} = (x+y\alpha)(x-y\alpha) = x^2 - y^2d$.
The norm mapping is a homomorphism from $\mathbb{F}^*_{q^2}$ to
$\mathbb{F}^*_{q}$ with $Im(N) = \mathbb{F}^*_q$. Thus, the kernel $Ker(N)$ is the subgroup of order $q+1$
in $\mathbb{F}_{q^2}^*$.

Now we make some remarks on squares in finite fields.
\begin{lem}\label{sq}
{\rm (1)} The element $-1$ is a square in $\mathbb{F}_q^*$ iff $q \equiv 1(4)$;\\
{\rm (2)} For any non-square $d$ in $\mathbb{F}_q^*$ the element $-d$ is a square in $\mathbb{F}_q^*$ iff $q \equiv 3(4)$.
\end{lem}

The following lemma can be used to test whether an element $\gamma = x+y\alpha \in \mathbb{F}_{q^2}^*$ is a square.

\begin{lem}[\cite{BEHW96}, Lemma 2]\label{sqq}
 An element $\gamma = x+y\alpha \in \mathbb{F}_{q^2}^*$ is a square iff $N(\gamma) = x^2-y^2d$ is a square in $\mathbb{F}_q^*$.
\end{lem}

Lemma \ref{sqqq} immediately follows from Lemma \ref{sqq}, Lemma \ref{sq} and the fact that $N(\alpha) = -d$.
\begin{lem}\label{sqqq}
The element $\alpha$ is a square in $\mathbb{F}_{q^2}^*$ iff $q \equiv 3(4)$.
\end{lem}

%

Denote by $\varphi$ the automorphism of $P(q^2)$ that sends $\gamma$ to $\gamma^q$.
Note that $\varphi$ fixes the elements from $\mathbb{F}_{q}$.

\medskip
Lemma \ref{phiSwaps} follows from Lemmas \ref{sqqq} and \ref{Frobenius}(2).
\begin{lem}\label{phiSwaps} The following statements hold.\\
{\rm(1)} If $q \equiv 1(4)$, then $\varphi$ interchanges only non-adjacent vertices.\\
{\rm(2)} If $q \equiv 3(4)$, then $\varphi$ interchanges only adjacent vertices.
\end{lem}
\bigskip

In view of Lemma \ref{phiSwaps}, and because the Paley graph $P(q^2)$ is isomorphic to its complement,
we can conclude that $P(q^2)$ satisfies the condition of Construction~\ref{C2} for every odd prime power $q$.
Thus if the parameters of $P(q^2)$ are $(4\mu +1,  2\mu, \mu -1, \mu)$, where $\mu = \frac{q^2-1}{4}$,
we can obtain the following three strictly Deza graphs from $P(q^2)$:
\\
(1) a strictly Deza graph with parameters $(4\mu +1,2\mu, \mu -1, \mu)$ obtained using dual Seidel switching in $P(q^2)$,
\\
(2) a strictly Deza graph with parameters $(8\mu +2,4\mu +1,4\mu,2\mu)$ obtained from Construction~\ref{C1},
\\
(3) a strictly Deza graph with parameters $(8\mu +2,4\mu +1,4\mu,2\mu)$ obtained from Construction~\ref{C2}.

\subsection{Symmetric conference matrices}

An $m\times m$ matrix $C$ with zero's on the diagonal, and $\pm 1$ elsewhere, is a {\em conference matrix} if $CC^\top=(m-1)I$.
If a conference matric $C$ is symmetric with constant row (and column) sum $r$, then $r=\pm\sqrt{m-1}$,
and $B=\frac{1}{2}(J_m-I_m-C)$ is the adjacency matrix of a strongly regular graph with parameter set
$${\cal P}(r) = \textstyle{ (\ r^2+1,\ \frac{1}{2}(r^2 - r),\ \frac{1}{4}(r - 1)^2 -1,\ \frac{1}{4}(r - 1)^2\ ) }.$$
Note that ${\cal P}(-r)$ is the complementary parameter set of ${\cal P}(r)$.
Symmetric conference matrices with constant row sum have been constructed by Seidel (see \cite{S}, Thm. 13.9).
If $q$ is an odd prime power and $r=\pm q$, then such a conference matrix can be obtained from the Paley graph of order $q^2$.
Let $B'$ be the adjacency matrix of $P(q^2)$, and put $S=J_{q^2}-I_{q^2}-2B'$ ($S$ is the so-called {\em Seidel matrix} of $P(q^2)$).
Define
$$C'=\left[\begin{array}{cc} 0 & {\bf 1}^\top\\ {\bf 1} & S \end{array}\right]$$
(${\bf 1}$ is the all-ones vector).
Then $C'$ is a symmetric conference matrix of order $m=q^2+1$.
However, $C'$ doesn't have constant row sum.

Next we shall make the row and column sum constant by multiplying some rows and the corresponding columns of $C'$ by $-1$.
This operation is called {\em Seidel switching}, and it is easily seen that Seidel switching doesn't change the conference matrix property.
To describe the required rows and columns, we use the notation and description of $P(q^2)$ given in the previous subsection.
If $q\equiv 3\mod 4$ we take the complement of the described Paley graph.
Then the involution $\varphi$ given in Lemma~\ref{phiSwaps} interchanges only non-adjacent vertices in all cases.
For $x\in\mathbb{F}_q$ define $V_x=\{x+y\alpha\ |\ y\in\mathbb{F}_q\}$.
Then the sets $V_x$ partition the vertex set of $P(q^2)$, and each class is a coclique.
Moreover, the partition is fixed by the involution $\varphi$.
Let $V$ be the union of $\frac{1}{2}(q-1)$ classes $V_x$.
Then $V$ induces a regular subgraph of $P(q^2)$ of degree $\frac{1}{4}(q-1)^2 -1$ with $\frac{1}{2}q(q-1)$ vertices.
Now make the matrix $C$ by Seidel switching in $C'$ with respect to the rows and columns that correspond with $V$.
Then $C$ is a regular symmetric conference matrix, and $B=\frac{1}{2}(J-I-C)$ is the adjacency matrix of a
strongly regular graph $\Gamma$ with parameter set ${\cal P}(q)$, and $\varphi$ remains an involution that
interchanges only nonadjacent vertices.
So $\Gamma$ satisfies the conditions for Construction~\ref{C1} and \ref{C2}.
For the complement of $\Gamma$ we found no involutions that interchanges only nonadjacent vertices.
Thus we find Deza graphs with parameters
$(q^2+1,\ \frac{1}{2}(q^2 - q),\ \frac{1}{4}(q - 1)^2,\ \frac{1}{4}(q - 1)^2-1\ )$ (by dual Seidel switching in $\Gamma$),
$(2q^2+2,\ q^2 - q+1),\ q^2 - q,\ \frac{1}{2}(q - 1)^2\ )$  (by Construction~\ref{C1} and \ref{C2} applied to $\Gamma$), and
$(2q^2+2,\ q^2 + q+1),\ q^2 + q,\ \frac{1}{2}(q + 1)^2\ )$  (by Construction~\ref{C1} applied to the complement of $\Gamma$).
If $q=3$, $\Gamma$ is the Petersen graph.
It has a unique involutive automorphism that interchanges only non-adjacent vertices.
This automorphism has 4 fixed and 6 moved vertices.
The Deza graph obtained from the Petersen graph with dual Seidel switching has diameter $3$.
However, in all other cases the obtained Deza graphs are strictly Deza.
The complement of the Petersen graph has no involutive automorphisms that interchanges only non-adjacent vertices.
We expect that this is the case for all strongly regular graphs with parameter set ${\cal P}(r)$ and $r<0$.

\subsection{The Hoffman-Singleton graph}
We consider Robertson's pentagons and pentagrams construction (see \cite{H03}) of Hoffman-Singleton graph, which is strongly regular
with parameters $(50, 7, 0, 1)$.

The 50 vertices of Hoffman-Singleton graph are grouped into 5 pentagons $P_0, \dots , P_4$ and 5 pentagrams $Q_0, \dots , Q_4$
labeled in such a way that the pentagrams are the complements of the pentagons; there are no edges between any two distinct pentagons,
nor between any two distinct pentagrams.
Edges between pentagon and pentagram vertices are defined by the rule:
each vertex $i\in \{0,1,2,3,4\}$ of a pentagon $P_j$, $j \in \{0,1,2,3,4\}$ is adjacent to the vertex $(i + jk)~mod~5$ of a
pentagram $Q_k$, for any $k \in \{0,1,2,3,4\}$.

Let $\varphi$ be the permutation of vertices of Hoffman-Singleton graph that fixes each vertex of $P_0$,  $Q_0$, and
interchanges $P_1$, $Q_1$ with $P_4$, $Q_4$ and $P_2$, $Q_2$ with $P_3$, $Q_3$.
The permutation $\varphi$ is a unique involutive automorphism of Hoffman-Singleton graph that interchanges only non-adjacent vertices.
The Deza graph obtained from Hoffman-Singleton graph with dual Seidel switching has diameter $3$.
However, each of Constructions~\ref{C1} and \ref{C2} produces a strictly Deza graph with parameters $(100, 15, 14, 2)$.
Construction~\ref{C1} applied to the complement gives a strictly Deza graph with parameters $(100,85,84,72)$.

\section{Divisible design graphs}

For a Deza graph $\Gamma$ with parameters $(n,k,b,a)$ and $a\neq b$,
we define two graphs $\Gamma_a$ and $\Gamma_b$ on the vertex set of $\Gamma$,
where two vertices $x$ and $y$ are adjacent in $\Gamma_a$ ($\Gamma_b$) if $x$ and $y$ have $a$ ($b$) common neighbors.
Clearly $\Gamma_a$ and $\Gamma_b$ are each others complement, and regular of degree $\alpha$ and $\beta$, respectively
($\Gamma_a$ and $\Gamma_b$ have been called the {\em children} of $\Gamma$).

If $\Gamma_a$ or $\Gamma_b$ is the disjoint union of complete graphs, then $\Gamma$ is called a {\em divisible design graph}
(DDG for short). DDGs are interesting structures on their own, and have been studied in \cite{HKM} and \cite{CH}.
If $A$ is the adjacency matrix of a DDG, then $A$ also satisfies the conditions for the incidence matrix of a divisible design,
which explains the name.

If $\Gamma$ is a Deza graph with $\alpha=1$ ($\beta=1$) then $\Gamma_a$ ($\Gamma_b$) consist of $n/2$ disjoint edges,
so $\Gamma$ is a DDG.
Such a DDG has been called {\em thin} (see \cite{CH}).
Thus Theorem~\ref{mainthm} also characterizes thin DDGs with $b=k-1$.
(for DDGs one uses $\lambda_1$ and $\lambda_2$ instead of $b$ and $a$).
The characterization of DDGs with $\lambda_1=k-1$ is also claimed in Theorem~4.11 of \cite{HKM}.
However, this claim is not correct, because it only uses Construction~\ref{C1} (which corresponds to Construction~4.10 in \cite{HKM}),
and Construction~\ref{C2} is not mentioned.
The proof in \cite{HKM} is based on the characterization of divisible designs with $k-\lambda_1=1$ from \cite{H91}.
But the authors of \cite{HKM} overlooked that isomorphic divisible designs may correspond to non-isomorphic DDGs.
Indeed, $A$ and $A'$ are adjacency matrices of isomorphic DDGs whenever there exist a permutation matrix $P$
such that $PAP^\top=A'$, but for $A$ and $A'$ to be incidence matrices of isomorphic divisible designs it is only required that
$PAQ=A'$ for permutation matrices $P$ and $Q$.
This is precisely what happens if one DDG can be obtained from the other by dual Seidel switching.
So Theorem~4.11 of \cite{HKM} can be repaired with basically the same proof by inserting Construction~\ref{C2} in the statement.

A graph is {\em walk-regular} if the number of closed walks of any given length at a vertex $x$ is independent of the choice of $x$.
From results in \cite{CH} it follows that the graphs made by Construction~\ref{C1} are walk regular, and those from Construction~\ref{C2}
are not.
It also follows that the matrices of Construction~\ref{C1} and \ref{C2} have different eigenvalues.
Because of this the discovered mistake in \cite{HKM} has some consequences for Table~1 in \cite{HKM} and \cite{CH},
and the second author likes to take the opportunity for giving the necessary corrections.
For the parameter sets $(18,9,8,4)$ and $(20,7,6,2)$, there can occur an additional eigenvalue $1$,
and the multiplicities are not determined by the parameters,
moreover for these parameter sets the "no" in column "notWR" (in Table~1 of \cite{CH}) should be replaced by a "yes".

\end{document}